\documentclass[a4paper,reqno,10pt]{amsart}
\usepackage{amssymb}
\usepackage{cite}
\usepackage{tikz}
\usepackage{amsfonts}
\usepackage{bigints}
\usepackage{amsmath}
\usepackage{amssymb}
\usepackage{systeme}
\usepackage{mathrsfs}
\usepackage[colorlinks=true]{hyperref}
\numberwithin{equation}{section}
\newtheorem{theorem}{Theorem}[section]
\newtheorem{prop}[theorem]{Proposition}
\newtheorem{lm}[theorem]{Lemma}
\newtheorem{cor}[theorem]{Corollary}

\newcommand{\R}{\mathbb{R}}

\def\ds{\displaystyle}

\def\dive{\mathrm{div}}
\newenvironment{preuve}{{\noindent {\bf Proof. }}}{\hfill {\rule{2.5mm}{2.5mm}}}

\makeatletter
\let\@msm@th@eqref\eqref
\renewcommand{\eqref}[1]{%
	\begingroup
	\leavevmode
	\color{red}%
	\hypersetup{linkbordercolor=[named]{red}}%
	\@msm@th@eqref{#1}%
	\endgroup
}
\makeatother

\author[M.~Amara]{Mustapha Amara}
\address{Department of Mathematics, Faculty of Science of Gab\`es, Research Laboratory Mathematics and Applications LR17ES11; Tunisia}
\email{\sl Mostafa.Amara@fsg.u-gabes.tn}

\title[Long-time behavior of global solutions of $(AQG)$ equations]
{Long-time behavior of global solutions of anisotropic quasi-geostrophic equations in Sobolev space}

\begin{document}
	\begin{abstract}
	We study the behavior at infinity in time of the global solution of the anisotropic quasi-geostrophic equation $\theta\in  C_b(\R^+,H^s( \R^2))$. We prove that this solution decays to zero as time goes to infinity in $L^p(\R^2)$, $p\in [2,+\infty],$ moreover, we prove also that $\ds\lim_{t\rightarrow +\infty}\|\theta(t)\|_{H^s}=0$.
	\end{abstract}
	
	

	\maketitle
	\tableofcontents

	
	\section{\bf Introduction and Main Results}
	We consider the following anisotropic quasi-geostrophic  equation given by \cite{YZ}:
	\begin{equation}\label{AQG}\tag{AQG}
		\begin{cases}
			\partial_t\theta+ u_\theta.\nabla\theta +\mu|\partial_1|^{2\alpha}\theta+\nu |\partial_2|^{2\beta}\theta=0,&\mbox{ in }\R^+\times\R^2\\
			\theta(x,0)=\theta^0(x),&\mbox{ in } \R^2
		\end{cases}
	\end{equation}
where $x\in\R^2$, $t>0$, $\alpha,\beta\in(0,1)$, $\mu,\nu>0$, $\theta=\theta(x,t)$ represents the potential temperature, and $u_\theta=(u_1,u_2)$ is the divergence free velocity which is determined by the Riesz transformation of $\theta$ in the following way:
	\begin{align*}
		&u_1:=-\mathcal{R}_2\theta=-\partial_{2}(-\Delta)^{-\frac{1}{2}}\theta,\\
		&u_2:=\mathcal{R}_1\theta=\partial_{1}(-\Delta)^{-\frac{1}{2}}\theta,
	\end{align*}

	The case when $\alpha=\beta$ and $\mu=\nu$, the system \eqref{AQG} back to the classical dissipation $(QG)$ equation:
		\begin{equation}\label{QG}\tag{QG}
		\begin{cases}
			\partial_t\theta+ u_\theta.\nabla\theta +\mu(-\Delta)^{2\alpha}\theta=0,&x\in\R^2, \ t>0\\
			\theta(x,0)=\theta^0(x).
		\end{cases}
	\end{equation}
	
	The equations \eqref{AQG} and \eqref{QG} are special cases of the general quasi-geostrophic approximations for atmospheric and oceanic fluid flow with small Rossby and Ekman numbers. We refer the reader to \cite{DC,CP,JP,CP1} where the authors explain the physical origin and the signification of the parameters of this equation.\\
	
	The existence of  global solution of the \eqref{AQG} system with $H^s(\R^2)$, $s\geq 2$, data is established by \cite{YZ} if $\alpha,\beta\in (0,1)$ satisfy
	\begin{equation}\label{1.1}
		\beta>\begin{cases}
			\frac{1}{2\alpha+1},&0<\alpha\leq \frac{1}{2}\\ \\
			
			\frac{1-\alpha}{2\alpha}&\frac{1}{2}<\alpha<1.
		\end{cases}
	\end{equation}
	In \cite{MJ}, Benameur and Amara proved that this system have a global solution when the initial data in $H^s(\R^2)$, $s\in (\max\{2-2\alpha,2-2\beta\},2)$, if $\alpha,\beta\in (1/2,1)$. We also refer to our result in \cite{MA}, where we show that the global solution of given by \cite{YZ}	and \cite{MJ} is bounded.
 \\
	
	The goal of this paper is the asymptotic study of this bounded global solution, this study of evolution equations has been a long time, see for example \cite{JB,JC}. In this paper, we will show that this global solution decay to zero as time goes to infinity in Lebesgue and Sobolev spaces.\\
	
	Before treating the decay rates of global solutions in the second section, we recall the results of global existence given in \cite{YZ,MJ,MA} which will be useful:
	\begin{theorem}[see \cite{YZ}]\label{Theorem1.1}
		Let $\alpha,\beta\in (0,1)$ and $s\geq 2$. If $\theta^0\in H^s(\R^2)$, then there exists a positive time $T>0$ depending on $\|\theta^0\|_{H^s}$ such that \eqref{AQG} admits a unique solution $\theta\in C([0,T],H^s(\R^2))$.
	\end{theorem}

Now we recall the following standard energy estimates.
\begin{prop}[see \cite{YZ}]
	Assume $\theta^0$ satisfies the assumptions stated in Theorem \ref{Theorem1.1} and let $\theta$ the corresponding solution. Then, for any $t> 0$,
	\begin{equation}
		\|\theta(t)\|_{L^2}^2+2\int_{0}^t\||\partial_{1}|^\alpha\theta\|_{L^2}^2d\tau+2\int_{0}^t\||\partial_{2}|^\beta\theta\|_{L^2}^2d\tau\leq \|\theta^0\|_{L^2}^2.
	\end{equation}
	Moreover, for any $p\in[2,+\infty]$
	\begin{equation}
		\|\theta(t)\|_{L^p}\leq \|\theta^0\|_{L^p}.
	\end{equation}
\end{prop}

\begin{theorem}[see \cite{MA}]\label{Theorem1.2}
	Let $\theta^0\in H^{s}(\R^2)$, $s\geq 2$, if $\alpha,\beta\in (0,1)$ satisfy \eqref{1.1}, then the system \eqref{AQG} admit a unique global bounded solution $\theta$ satisfy
	$$\theta\in C_b(\R^+,H^s(\R^2)),\quad|\partial_{1}|^\alpha\theta,|\partial_{2}|^\beta\theta\in L^2(\R^+,H^s(\R^2)).$$
	Moreover, for any $t>0$ we have 
	\begin{equation}
		\|\theta(t)\|_{H^s}^2+\int_0^t\||\partial_{1}|^\alpha\theta\|_{H^s}^2d\tau+\int_0^t\||\partial_{2}|^\beta\theta\|_{H^s}^2d\tau\leq C(\theta^0),
	\end{equation}
where $C(\theta^0)$ is a constant depending on the initial data $\theta^0$.
\end{theorem}
\begin{theorem}[see \cite{MA}]\label{theorem1.4}
	Let $\alpha,\beta\in (1/2,1)$ and $\theta^0\in H^{s}(\R^2)$, $s>\max\{2-2\alpha,2-2\beta\}$. Then the system \eqref{AQG} admit a unique global bounded solution $\theta$ satisfy
	$$\theta\in C_b(\R^+,H^s(\R^2))\cap C_b((0,+\infty),H^2(\R^2)),\quad|\partial_{1}|^\alpha\theta,|\partial_{2}|^\beta\theta\in L^2(\R^+,H^s(\R^2)).$$
\end{theorem}
Now we are ready to state our first result:
\begin{theorem}\label{theorem3}
	Let $\theta\in (\R^+,H^s(\R^2))$, $s\geq 2$, be a global solution of \eqref{AQG}, then, for any $p\in[ 2,+\infty)$ we have 
	\begin{equation}
		\lim_{t\rightarrow +\infty}\|\theta(t)\|_{L^p}=0.
	\end{equation}
Moreover, if $\alpha,\beta$ satisfies \eqref{1.1}, therefore
\begin{equation}
	\lim_{t\rightarrow +\infty}\|\theta(t)\|_{L^{\infty}}=0.
\end{equation}
\end{theorem}

Moving to announce our second result.
\begin{theorem}\label{theorem4}
		Let $\theta\in C_b(\R^+,H^s(\R^2))$, $s\geq2$, be the global solution of \eqref{AQG} and $\alpha,\beta$ satisfy \eqref{1.1}. Then
	\begin{equation}
		\lim_{t\rightarrow +\infty}\|\theta(t)\|_{H^s}=0.
	\end{equation}
\end{theorem}
The next corollary show that if our result hold true in the Sobolev space $H^s(\R^2)$, when $s> \max\{2-2\alpha,2-2\beta\}$.
\begin{cor}\label{Corollary2.5}
	Let	$\alpha,\beta\in (1/2,1)$ and  $\theta^0\in H^s(\R^2)$, $s\in (\max\{2-2\alpha,2-2\beta\},2) $ be the global solution of \eqref{AQG}. Then
	\begin{equation}
		\lim_{t\rightarrow +\infty}\|\theta(t)\|_{H^s}=0.
	\end{equation}
\end{cor}
\section{\bf Notation and Preliminary Results}
	In this short section, we collect some notations and definitions that will be used later, and we	give some technical lemmas.
	\begin{itemize}
		\item[$\bullet$] The Fourier transformation in $\R^2$		
		\begin{equation}
			\mathcal{F}(f)(\xi)=\widehat{f}(\xi)=\int_{\R^2}e^{-ix.\xi}f(x)dx,\quad \xi\in \R^2.
		\end{equation}
		The inverse Fourier formula is
		\begin{equation}
			\mathcal{F}^{-1}(f)(x)=(2\pi)^{-2}\widehat{f}(\xi)=\int_{\R^2}e^{i\xi.x}f(\xi)d\xi,\quad x\in \R^2.
		\end{equation}
		\item[$\bullet$] The convolution product of a suitable pair of function $f$ and $g$ on $\R^2$ 	is given by
		\begin{equation}
			f\ast g(x)=\int_{\R^2} f(x-y)g(y)dy.
		\end{equation}
		\item[$\bullet$] If $f=(f_1,f_2)$ and $g=(g_1,g_2)$ are two vector fields, we set
		$$f\otimes g:= (g_1f,g_2f),$$
		and 
		$$\dive (f\otimes g):= (\dive(g_1f),\dive(g_2f)).$$
		\item[$\bullet$] For $s\in \R$:
		\begin{itemize}
			\item[$\ast$] The Sobolev space: 
			\begin{align*}
				&H^{s}(\R^2):=\left\{f\in \mathcal{S}'(\R^2); (1+|\xi|^2)^{s/2}\widehat{f}\in L^2(\R^2)\right\},
			\end{align*} denotes the usual inhomogeneous Sobolev space on $\R^2$, with the norm 
			\begin{align*}
				&	\|f\|_{H^s}=\left(\int_{\R^2}(1+|\xi|^2)^s|\widehat{f}(\xi)|^2d\xi\right)^{\frac{1}{2}},
			\end{align*}
		and the scalar product 
		\begin{align*}
			&	\left(f,g\right)_{H^s}=\int_{\R^2}(1+|\xi|^2)^s\widehat{f}(\xi)\overline{\widehat{g}(\xi)}d\xi.
		\end{align*}
		\item[$\ast$] The Sobolev space: 
		\begin{align*}
			 \dot{H}^{s}(\R^2):=\left\{f\in \mathcal{S}'(\R^2);\widehat{f}\in L^1(\R^2)\mbox{ and }|\xi|^s\widehat{f}\in L^2(\R^2)\right\},
		\end{align*} denotes the usual  homogeneous Sobolev space on $\R^2$, with the norm 
		\begin{align*}
			\|f\|_{\dot{H}^s}=\left(\int_{\R^2}|\xi|^{2s}|\widehat{f}(\xi)|^2d\xi\right)^{\frac{1}{2}},
		\end{align*}
		and the scalar product 
	\begin{align*}
		&	\left(f,g\right)_{\dot{H}^s}=\int_{\R^2}|\xi|^{2s}\widehat{f}(\xi)\overline{\widehat{g}(\xi)}d\xi.
	\end{align*}
		\end{itemize}
	\end{itemize}

	We recall a fundamental lemma concerning the Sobolev spaces
		We recall a fundamental lemma concerning the Sobolev spaces
	\begin{lm}[see \cite{MA}]\label{Lemma}
		For $s,s_1,s_2\in\R$ and $z\in[0,1]$, the following anisotropic interpolation inequalities hold true for $i=1,2$:
		\begin{align}
			\label{ing1}	&\||\partial_{i}|^{zs_1+(1-z)s_2}f\|_{H^{s}}\leq \||\partial_{i}|^{s_1}f\|_{H^{s}}^z\||\partial_{i}|^{s_2}f\|_{H^{s}}^{1-z},\\
			\label{ing2}	&\||\partial_{i}|^{zs_1+(1-z)s_2}f\|_{\dot{H}^{s}}\leq \||\partial_{i}|^{s_1}f\|_{\dot{H}^{s_1}}^z\||\partial_{i}|^{s_2}f\|_{\dot{H}^{s_2}}^{1-z}.
		\end{align}
	\end{lm}
	\begin{lm}[see \cite{BH}]\label{Lemma1}
		Let $s_1$, $s_2$ be two real numbers such that $s_1<1$ and $s_1+s_2>0$. Then, there exists a positive constant $C=C(s_1,s_2)$ such that for all $f,g\in \dot{H}^{s_1}(\R^2)\bigcap \dot{H}^{s_2}(\R^2)$; 
		\begin{equation}
			\|fg\|_{\dot{H}^{s_1+s_2-1}}\leq C(s_1,s_2) \left(\|f\|_{\dot{H}^{s_1}}\|g\|_{\dot{H}^{s_2}}+\|f\|_{\dot{H}^{s_2}}\|g\|_{\dot{H}^{s_1}}\right).
		\end{equation}
		Moreover, if $s_2<1$, there exists a positive constant $C'=C'(s_1,s_2)$ such that for all $f\in \dot{H}^{s_1}(\R^2)$ and $g\in \dot{H}^{s_2}(\R^2)$; 
		\begin{equation}
			\|fg\|_{\dot{H}^{s_1+s_2-1}}\leq C'(s_1,s_2)\|f\|_{\dot{H}^{s_1}}\|g\|_{\dot{H}^{s_2}}
		\end{equation}
	\end{lm}
	\begin{lm}[see \cite{JN}]\label{Lemma2}
		For any $p\in (1,+\infty)$, there is a constant $C(p)>0$ such that
		\begin{equation}
			\|\mathcal{R}^\perp\theta\|_{L^p}\leq C(p) \|\theta\|_{L^p}.
		\end{equation}
	\end{lm}
	\par We recall the following important commutator and product estimates:
	\begin{lm}[see \cite{MJ}]\label{Lemma4}
		Suppose that $s>1$. If $f,g\in \mathcal{S}$ then for any $\alpha\in(0,1)$
		\begin{align*}
			\||\nabla|^s(f g)-f|\nabla|^s g\|_{L^2}\leq C\left(\||\nabla|^{s+\alpha} f\|_{L^{2}}\||\nabla|^{1-\alpha}g\|_{L^{2}}+\||\nabla|^{s-1+\alpha} g\|_{L^{2}}\||\nabla|^{2-\alpha}f\|_{L^{2}}\right)
		\end{align*}
		and
		\begin{align*}
			\||\nabla|^s(f g)\|_{L^2}\leq C\left(\||\nabla|^{s+\alpha} f\|_{L^{2}}\||\nabla|^{1-\alpha}g\|_{L^{2}}+\||\nabla|^{s+\alpha} g\|_{L^{2}}\||\nabla|^{1-\alpha}f\|_{L^{2}}\right).
		\end{align*}
	\end{lm}
\begin{lm}\label{Lemma5}
	Let $r>0$, then, there is a constant $C(r)$ such that
	\begin{equation}
	\sup_{\xi\in \R^2}\left(	|\xi_i|^{r} \left|\widehat{f g}(\xi)\right|\right) \leq C(r)\left(\|f\|_{L^2}\||\partial_i|^rg\|_{L^2}+\|g\|_{L^2}\||\partial_i|^rf\|_{L^2}\right),\quad\forall i\in\{1,2\}.
	\end{equation}
\end{lm}

\begin{preuve} We have for any $\xi\in \R^2$ and $i\in \{1,2\}$
		\begin{align*}
		|\xi_i|^r \left|\widehat{fg}(\xi)\right|&\leq \int_{\R^2} |\xi_i|^r |\widehat{f}(\xi-\eta)|\widehat{g}(\eta)|d\eta\\
		& \leq C(r)  \int_{\R^2} \left(|\xi_i-\eta_i|^r+|\eta_i|^r\right)  |\widehat{f}(\xi-\eta)|\widehat{g}(\eta)|d\eta\\
		& \leq C(r)\int_{\R^2} \left(|\widehat{f_1}(\xi-\eta)|\widehat{g}(\eta)+|\widehat{f}(\xi-\eta)|\widehat{g_1}(\eta)|\right)d\eta\\
		& \leq C(r) \left(\|\widehat{f_1g}\|_{L^\infty}+\|\widehat{fg_1}\|_{L^\infty}\right),
	\end{align*}
where 
$$\widehat{f_1}(\xi)=|\xi_i|^r\left|\widehat{f}(\xi)\right|\mbox{ and } \widehat{g_1}(\xi)=|\xi_i|^r\left|\widehat{g}(\xi)\right|.$$
Therefore 
\begin{align*}
\sup_{\xi\in \R^2}\left(	|\xi_i|^r \left|\widehat{fg}(\xi)\right|\right)&\leq  C(r) \left(\|f_1g\|_{L^1}+\|fg_1\|_{L^1}\right)\\
&\leq  C(r) \left(\|f_1\|_{L^2}\|g\|_{L^2}+\|f\|_{L^2}\|g_1\|_{L^2}\right)\\
&\leq  C(r) \left(\|g\|_{L^2}\||\partial_i|^rf\|_{L^2}+\|f\|_{L^2}\||\partial_i|^rg\|_{L^2}\right).
\end{align*}
\end{preuve}
	\section{\bf Proof of Theorem \ref{theorem3}}
	The decay study in $L^2$ is inspired by \cite{JR}, precisely we cut $\theta$ into high and low frequencies. For that, let $\delta$ a strictly positive real number strictly less than 1, we define the operators
		$A_\delta(D)$ and $B_\delta(D)$, respectively, by the following:
		\begin{equation}
			\begin{array}{l}
				A_\delta(D)f=\mathcal{F}^{-1}(\chi_{\mathscr{P}(0,\delta)}\mathcal{F}(f))\mbox{ and }
				B_\delta(D)f=\mathcal{F}^{-1}((1-\chi_{\mathscr{P}(0,\delta)})\mathcal{F}(f)),
			\end{array}
		\end{equation}
		where $$\mathscr{P}(0,\delta)=\left\{\xi=(\xi_1,\xi_2)\in \R^2,\ \max\{|\xi_1|,|\xi_2|\}\leq \delta\right\}$$ and 
		$$\mathscr{P}(0,\delta)^c=\left\{\xi=(\xi_1,\xi_2)\in \R^2,\ \max\{|\xi_1|,|\xi_2|\}> \delta\right\}.$$
		We define $w_\delta=A_\delta(D)\theta$ and $v_\delta=B_\delta(D)\theta$; $\mathcal{F}(\theta)=\mathcal{F}(w_\delta)+\mathcal{F}(v_\delta)$. Then,
		\begin{equation}\label{5.2}
			\partial_tw_\delta+|\partial_1|^{2\alpha}w_\delta+|\partial_2|^{2\beta}w_\delta+A_\delta(D)( u_\theta.\nabla\theta)=0
		\end{equation}
		and
		\begin{equation}
			\partial_tv_\delta+|\partial_1|^{2\alpha}v_\delta+|\partial_2|^{2\beta}v_\delta+B_\delta(D)( u_\theta.\nabla\theta)=0
		\end{equation}
		Taking the scalar product of \eqref{5.2} equation with $w_\delta$, we get
		$$	\frac{1}{2}	\frac{d}{dt}\|w_\delta(t)\|_{L^2}^2+\||\partial_{1}|^\alpha w_\delta(t)\|_{L^2}^2+\||\partial_{2}|^\beta w_\delta(t)\|_{L^2}^2\leq \left|\left(A_\delta(D)(u_\theta.\nabla\theta),w_\delta\right)_{L^2}\right|.$$ We have
		\begin{align*}
			\left|\left(A_\delta(D)(u_\theta.\nabla\theta),w_\delta\right)_{L^2}\right|
			&\leq \int_{\mathscr{P}(0,\delta)}|\xi||\widehat{u_\theta  \theta}(\xi)||\widehat{w_\delta}(\xi)|d\xi\\
			&\leq \int_{\mathscr{F}_1}|\xi||\widehat{u_\theta  \theta}(\xi)||\widehat{w_\delta}(\xi)|d\xi+\int_{\mathscr{F}_2}|\xi||\widehat{u_\theta  \theta}(\xi)||\widehat{w_\delta}(\xi)|d\xi\\
			&\leq 2\int_{\mathscr{F}_1}|\xi_1||\widehat{u_\theta  \theta}(\xi)||\widehat{w_\delta}(\xi)|d\xi+2\int_{\mathscr{F}_2}|\xi_2||\widehat{u_\theta  \theta}(\xi)||\widehat{w_\delta}(\xi)|d\xi\\
			&\leq 2\int_{\mathscr{F}_1}|\xi_1|^{1-2\alpha}|\xi_1|^\alpha|\widehat{u_\theta  \theta}(\xi)||\xi_1|^\alpha|\widehat{w_\delta}(\xi)|d\xi+2\int_{\mathscr{F}_2}|\xi_2|^{1-2\beta}|\xi_2|^\beta|\widehat{u_\theta \theta}(\xi)||\xi_2|^\beta|\widehat{w_\delta}(\xi)|d\xi,
		\end{align*}
	where $\mathscr{F}_1:=\mathscr{P}(0,\delta)\bigcap \left\{\xi\in \R^2;\ |\xi_2|\leq |\xi_1|\right\}$ and $\mathscr{F}_2:=\mathscr{P}(0,\delta)\bigcap \left\{\xi\in \R^2;\ |\xi_1|< |\xi_2|\right\}$.\\
		
	By the Lemma \ref{Lemma5} we have
		\begin{align*}
			|\xi_1|^\alpha |\widehat{u_\theta  \theta}(\xi)|\leq C\|\theta\|_{L^2}\||\partial_{1}|^\alpha\theta\|_{L^2}\leq C\|\theta^0\|_{L^2}\||\partial_{1}|^\alpha\theta\|_{L^2}, 
		\end{align*}
		and
		\begin{align*}
			|\xi_1|^\beta |\widehat{u_\theta  \theta}(\xi)|\leq C\|\theta\|_{L^2}\||\partial_{2}|^\beta\theta\|_{L^2}\leq C\|\theta^0\|_{L^2}\||\partial_{2}|^\beta\theta\|_{L^2}.
		\end{align*}
		So we get 
		\begin{align*}
			\left|\left(A_\delta(D)(u_\theta.\nabla\theta),w_\delta\right)_{L^2}\right|
			&\leq C\|\theta^0\|_{L^2}\||\partial_{1}|^\alpha\theta\|_{L^2}\int_{\mathscr{F}_1}|\xi_1|^{1-2\alpha}|\widehat{|\partial_{1}|^\alpha w_\delta}(\xi)|d\xi\\
			&\hskip1cm+C\|\theta^0\|_{L^2}\||\partial_{2}|^\beta\theta\|_{L^2}\int_{\mathscr{F}_2}|\xi_2|^{1-2\beta}|\widehat{|\partial_{2}|^\beta w_\delta}(\xi)|d\xi
		\end{align*}
		Moreover, we have 
		\begin{align*}
			\int_{\mathscr{F}_1}|\xi_1|^{1-2\alpha}|\widehat{|\partial_{1}|^\alpha w_\delta}(\xi)|d\xi&\leq \begin{cases}
				\ds\delta^{1-2\alpha} \int_{\mathscr{B}(0,2\delta)}|\widehat{|\partial_{1}|^\alpha w_\delta}(\xi)|d\xi,&\mbox{ if } \alpha\in (0,1/2],\\
				\\
				\ds 2\int_{\mathscr{B}(0,2\delta)}|\xi|^{1-2\alpha}|\widehat{|\partial_{1}|^\alpha w_\delta}(\xi)|d\xi,&\mbox{ if } \alpha\in (1/2,1),
			\end{cases}
		\end{align*}
	where $\mathscr{B}(0,2\delta)=\left\{\xi\in \R^2,\ |\xi|<2\delta\right\}.$\\
	
	 We start with the case when $\alpha\in (0,1/2]$, and by Cauchy-Schwartz inequality we have 
	 \begin{align*}
	\delta^{1-2\alpha} 	\int_{\mathscr{B}(0,2\delta)}|\widehat{|\partial_{1}|^\alpha w_\delta}(\xi)|d\xi&\leq \delta^{1-2\alpha}  \left(\int_{\mathscr{B}(0,2\delta)}d\xi\right)^{1/2}\||\partial_{1}|^\alpha w_\delta\|_{L^2}\leq C\delta^{2-2\alpha}  \||\partial_{1}|^\alpha \theta\|_{L^2}.
	 \end{align*}
 Now move for the case if $\alpha\in (1/2,1)$, we have also by Cauchy-Schwartz inequality
		\begin{align*}
			\int_{\mathscr{B}(0,2\delta)}|\xi|^{1-2\alpha}|\widehat{|\partial_{1}|^\alpha w_\delta}(\xi)|d\xi&\leq \left(\int_{\mathscr{B}(0,2\delta)}|\xi|^{2-4\alpha}d\xi\right)^{1/2}\||\partial_{1}|^\alpha w_\delta\|_{L^2}\\
			&\leq C \left(\int_{0}^{2\delta}r^{3-4\alpha}dr\right)^{1/2}\||\partial_{1}|^\alpha \theta\|_{L^2}\leq C\delta^{2-2\alpha}\||\partial_{1}|^\alpha \theta\|_{L^2}.
		\end{align*}
	Therefore, for any $\alpha\in(0,1)$ we have 
	\begin{align*}
			\int_{\mathscr{F}_1}|\xi_1|^{1-2\alpha}|\widehat{|\partial_{1}|^\alpha w_\delta}(\xi)|d\xi\leq C \delta^{2-2\alpha}\||\partial_{1}|^\alpha \theta\|_{L^2}.
	\end{align*}
Same thing for the other term, for any $\beta\in(0,1)$, we have 
		\begin{align*}
			\int_{\mathscr{F}_2}|\xi_2|^{1-2\beta}|\widehat{|\partial_{2}|^\beta w_\delta}(\xi)|d\xi&\leq C(\beta)\delta^{2-2\beta}\||\partial_{2}|^\beta \theta\|_{L^2}.
		\end{align*}
		Finally, we get
		\begin{align}
			\frac{1}{2}	\frac{d}{dt}\|w_\delta(t)\|_{L^2}^2\leq  C \left(\delta^{2-2\alpha}+\delta^{2-2\beta}\right)\|\theta^0\|_{L^2}\left(\||\partial_{1}|^\alpha\theta\|_{L^2}^2+\||\partial_{2}|^\beta\theta\|_{L^2}^2\right).
		\end{align}
		Integer in $[0,t]$, $t>0$, we get
		\begin{align*}
			\|w_\delta(t)\|_{L^2}^2&\leq 	\|w_\delta^0\|_{L^2}^2+C \left(\delta^{2-2\alpha}+\delta^{2-2\beta}\right)\|\theta^0\|_{L^2} \int_{0}^{t}\left(\||\partial_{1}|^\alpha\theta\|_{L^2}^2+\||\partial_{2}|^\beta\theta\|_{L^2}^2\right)d\tau\\
			&\leq \|w_\delta^0\|_{L^2}^2+ C \left(\delta^{2-2\alpha}+\delta^{2-2\beta}\right)\|\theta^0\|_{L^2}^3\underset{\delta\rightarrow0^+}{\longrightarrow}0.
		\end{align*}
		Which implies
		\begin{equation*}
			\lim\limits_{\delta\rightarrow 0^+}\|w_\delta(t)\|_{L^\infty(\R_+,L^2(\R^2))}=0.
		\end{equation*}
		Let $\varepsilon>0$, then there exists $\delta_\varepsilon>0$ such that
	\begin{equation}\label{5.5}
		\|w_{\delta_\varepsilon}\|_{L^\infty(\R_+,L^2(\R^2))}<\frac{\varepsilon}{2}.
	\end{equation}
		On the other hand, we have 
		\begin{align*}
			\|v_\delta\|_{L^2}^2=\int_{\mathscr{P}(0,\delta)^c} |\widehat{\theta}(\xi)|^2d\xi
			&= \int_{\mathscr{F}_3}|\widehat{\theta}(\xi)|^2d\xi+\int_{\mathscr{F}_4}|\widehat{\theta}(\xi)|^2d\xi,
		\end{align*}
		where $\mathscr{F}_3:=\mathscr{P}(0,\delta)^c\bigcap \left\{\xi\in \R^2;\ |\xi_1|\leq |\xi_2|\right\}$ and $\mathscr{F}_2:=\mathscr{P}(0,\delta)^c\bigcap \left\{\xi\in \R^2;\ |\xi_2|< |\xi_1|\right\}$.
	 Therefore
	\begin{align*}
		\int_{\mathscr{F}_1}|\widehat{\theta}(\xi)|^2d\xi&=\int_{\mathscr{F}_1}|\xi_2|^{-2\beta}|\xi_2|^{2\beta}|\widehat{\theta}(\xi)|^2d\xi\leq \delta^{-2\beta}\||\partial_{2}|^\beta\theta\|_{L^2}^2.
	\end{align*}
The something for other term 
	\begin{align*}
	\int_{\mathscr{F}_2}|\widehat{\theta}(\xi)|^2d\xi&=\int_{\mathscr{F}_2}|\xi_1|^{-2\alpha}|\xi_1|^{2\alpha}|\widehat{\theta}(\xi)|^2d\xi\leq \delta^{-2\alpha}\||\partial_{1}|^\alpha\theta\|_{L^2}^2.
\end{align*}
Therefore 
	\begin{align*}
	\|v_\delta\|_{L^2}^2&\leq \frac{1}{\delta^{2\alpha}}\||\partial_{1}|^\alpha\theta\|_{L^2}^2+\frac{1}{\delta^{2\beta}}\||\partial_{2}|^\beta\theta\|_{L^2}^2.
\end{align*}
Now, we consider $S(\delta_\varepsilon)=\{t\geq 0; \|v_{\delta_\varepsilon}\|_{L^2}>\frac{\varepsilon}{2}\}$ the sub-part of $\R^+$, then $S(\delta_\varepsilon)$ is of finite measure, in fact
		\begin{align*}
			\left(\frac{\varepsilon}{2}\right)^2\lambda(S(\delta_\varepsilon))&\leq \int_{S(\delta_\varepsilon)}\|v_{\delta_\varepsilon}\|_{L^2}^2dt\\
			&\leq \int_{0}^{+\infty}\|v_{\delta_\varepsilon}\|_{L^2}^2dt\\
			&\leq\int_0^{+\infty}\left(\frac{1}{\delta_\varepsilon^{2\alpha}}\||\partial_{1}|^\alpha\theta\|_{L^2}^2+\frac{1}{\delta_\varepsilon^{2\beta}}\||\partial_{2}|^\beta\theta\|_{L^2}^2\right)dt\leq \left(\frac{1}{\delta_\varepsilon^{2\alpha}}+\frac{1}{\delta_\varepsilon^{2\beta}}\right)\|\theta^0\|_{L^2}^2.
		\end{align*}
		We pose
		$$T_{\varepsilon}=\left(\frac{2}{\varepsilon}\right)^2\left(\frac{1}{\delta_\varepsilon^{2\alpha}}+\frac{1}{\delta_\varepsilon^{2\beta}}\right)\|\theta^0\|_{L^2}^2<+\infty,$$
		then $\lambda(S_\varepsilon(\delta_0))\leq T_{\varepsilon}$. So there exists $t_\varepsilon\in [0,T_\varepsilon+1]\setminus S(\delta_\varepsilon)$ such
		\begin{equation}\label{5.6}
			\|v_{\delta_0}(t_\varepsilon)\|_{L^2}\leq \frac{\varepsilon}{2}.
		\end{equation}
		By the equation \eqref{5.5} and \eqref{5.6}, we get
		\begin{equation*}
			\|\theta(t_\varepsilon)\|_{L^2}\leq \varepsilon.
		\end{equation*}
		The fact that $\left(t\mapsto\theta(t)\right)$  is a decreasing function in $L^2$. Therefore
		 \begin{equation}
		 	\lim\limits_{t\rightarrow+\infty}\|\theta(t)\|_{L^2}=0.
		 \end{equation}
	Passing now to show that the result is true on all Lebesgue spaces $L^p(\R^2)$, $p\in(2,+\infty)$. So, taking $p\in (2,+\infty)$ and we have
		\begin{align*}
			\|\theta(t)\|_{L^p}^p&= \int_{\R^2}\left|\theta(x,t)\right|^p dx\\
			&= \int_{\R^2}\left|\theta(x,t)\right|^{p-2}\left|\theta(x,t)\right|^{2}  dx\\
			&\leq \|\theta(t)\|_{L^\infty}^{p-2} \|\theta(t)\|_{L^2}^2\\
			&\leq \|\theta^0\|_{L^\infty}^{p-2} \|\theta(t)\|_{L^2}\underset{t\rightarrow+\infty}{\longrightarrow}0.
		\end{align*}
Therefore, for any $p\in[2,+\infty)$
\begin{equation*}
	\lim_{t\rightarrow +\infty}\|\theta(t)\|_{L^p}=0.
\end{equation*}
Finally, the fact that $H^{s'}(\R^2)\hookrightarrow L^\infty(\R^2)$, $1<s'<s$, we have, for any $t\geq 0$
\begin{align*}
	\|\theta(t)\|_{L^\infty}&\leq C\|\theta(t)\|_{H^{s'}}\leq C\|\theta(t)\|_{L^2}^{1-\frac{s'}{s}}\|\theta(t)\|_{H^s}^{\frac{s'}{s}},
\end{align*}
where we use Lemma \ref{Lemma} for $0<s'<s$.\\

If $\alpha$ and $\beta$ satisfy \eqref{1.1} then by the Theorem \ref{Theorem1.2}; we have $\|\theta(t)\|_{H^s}\geq C(\theta^0)$ and according to the above; we get
\begin{equation}
\|\theta(t)\|_{L^\infty}\leq C(\theta^0)\|\theta(t)\|_{L^2}^{1-\frac{s'}{s}}\underset{t\rightarrow+\infty}{\longrightarrow}0.
\end{equation}
Therefore 
\begin{equation*}
	\lim_{t\rightarrow +\infty}\|\theta(t)\|_{L^\infty}=0,
\end{equation*}
which implies the desired result.
\hfill$\blacksquare$

	\section{\bf Proof of Theorem \ref{theorem4}}
The proof is done in two steps; the first step is to show that for any small constant, there exists a positive time such that the norm of our solution in this time is lower than this constant. In the second step, we consider our system with initial condition in this time and show that the solution is decreasing. This result shows that our solution tends to zero in the neighborhood of infinity.
\subsection{Step 1:}
Using Lemma \ref{Lemma}, $\alpha<s<s+\alpha$ and $\beta<s<s+\beta$, we get
	\begin{align*}
		\|\theta(t)\|_{\dot{H}^s}^2&\leq C\left( \||\partial_{1}|^{s}\theta\|_{L^2}^2+\||\partial_{2}|^{s}\theta\|_{L^2}^2\right)\\
		&\leq C\left( \||\partial_{1}|^{\alpha}\theta(t)\|_{H^s}^2+\||\partial_{2}|^{\beta}\theta(t)\|_{H^s}^2\right).
	\end{align*}
Therefore, by the Theorem \ref{Theorem1.1},
\begin{align*}
	\int_0^{+\infty}\|\theta(t)\|_{\dot{H}^s}^2dt&\leq C \int_0^{+\infty}\left( \||\partial_{1}|^{\alpha}\theta(t)\|_{H^s}^2+\||\partial_{2}|^{\beta}\theta(t)\|_{H^s}^2\right)dt \leq C(\theta^0)<+\infty.
\end{align*}
Let $\varepsilon>0$ and  considering the set  $S(\varepsilon)=\{t\geq 0; \|\theta(t)\|_{\dot{H}^s}>\varepsilon\}$, then
	\begin{align*}
		\varepsilon^2\lambda(S(\varepsilon))&\leq \int_{S(\varepsilon)}\|{\theta(t)}\|_{\dot{H}^s}^2dt\\
		&\leq \int_{0}^{+\infty}\|\theta(t)\|_{\dot{H}^s}^2dt\\
		&\leq C(\theta^0).
	\end{align*}
	We pose
	$$T_{\varepsilon}=\frac{C(\theta^0)}{\varepsilon^2}<+\infty,$$
	then $\lambda(S(\varepsilon))\leq T_{\varepsilon}$. So there exists $t_\varepsilon\in [0,T_\varepsilon+1]\setminus S(\varepsilon)$ such that
	\begin{equation}\label{6.1}
		\|\theta(t_\varepsilon)\|_{\dot{H}^s}\leq \varepsilon.
	\end{equation}
\subsection{Step 2:} Let $\varepsilon>0$, then, by the first step, there exists $t_\varepsilon\geq 0$ such that
\begin{equation*}
	\|\theta(t_\varepsilon)\|_{\dot{H}^s}\leq \varepsilon.
\end{equation*}
Considering now the follow system
\begin{equation}\label{AQG1}\tag{$AQG_1$}
	\begin{cases}
		\partial_t\gamma+ u_\gamma.\nabla\gamma +|\partial_1|^{2\alpha}\gamma+ |\partial_2|^{2\beta}\theta=0,\\
		\gamma(0)=\theta(t_\varepsilon).
	\end{cases}
\end{equation}
By Theorem \ref{Theorem1.1}, there exist a unique global solution $\gamma\in C_b(\R^+,H^s)$. The uniqueness of the solution show that
$$\gamma(t)=\theta(t+t_\varepsilon),\quad\forall t\geq 0.$$
Moreover, we have by applying $1+|\nabla|^s$ to \eqref{AQG1} and taking $L^2$ inner product with $(1+|\nabla|^s)\gamma$ 
\begin{align*}
	\frac{1}{2}\frac{d}{dt}\|\gamma\|_{H^s}^2+\||\partial_{1}|^\alpha\gamma\|_{H^s}^2+\||\partial_{2}|^\beta\gamma\|_{H^s}^2\leq \||\nabla|^s(u_\gamma.\nabla\gamma)-u_\gamma|\nabla|^s\nabla\gamma\|_{\dot{H}^s}\|\gamma\|_{\dot{H}^s},
\end{align*}
where we use $\dive(u_\gamma)=0$. Moreover by lemma \ref{Lemma4} we have 
\begin{align*}
	\frac{1}{2}\frac{d}{dt}\|\gamma(t)\|_{H^s}^2+\||\partial_{1}|^\alpha\gamma\|_{H^s}^2+\||\partial_{2}|^\beta\gamma\|_{H^s}^2\leq C\left(\||\partial_{1}|^\alpha\gamma\|_{H^s}^2+\||\partial_{2}|^\beta\gamma\|_{H^s}^2\right)\|\gamma\|_{\dot{H}^s}.
\end{align*}
Therefore, for $T_*=\sup\left\{T\geq 0;\ \ds\sup_{0\leq t\leq T}\|\gamma(t)\|_{H^s}\leq 2\varepsilon\right\}$. Take $t\in [0,T_*] $ Then we have
\begin{align*}
	\frac{1}{2}\frac{d}{dt}\|\gamma(t)\|_{H^s}^2+\||\partial_{1}|^\alpha\gamma\|_{H^s}^2+\||\partial_{2}|^\beta\gamma\|_{H^s}^2&\leq  2C\left(\||\partial_{1}|^\alpha\gamma\|_{H^s}^2+\||\partial_{2}|^\beta\gamma\|_{H^s}^2\right)\varepsilon.
\end{align*}
Taking $\varepsilon<\frac{1}{4C}$, and integrating over $(0, t)$, $t\in [0,T_*)$, we get
\begin{align*}
	\|\gamma(t)\|_{H^s}^2+\int_{0}^t\||\partial_{1}|^\alpha\gamma\|_{H^s}^2d\tau+\int_{0}^t\||\partial_{2}|^\beta\gamma\|_{H^s}^2d\tau&\leq \|\gamma(0)\|^2_{H^s}<2\varepsilon.
\end{align*}
Therefore $T_*=+\infty$ and we have for any $t>t_\varepsilon$
\begin{align*}
	\|\theta(t)\|_{\dot{H}^s}\leq \varepsilon,
\end{align*}
which implies that
$$\lim_{t\rightarrow +\infty} \|\theta(t)\|_{\dot{H}^s}=0.$$
By Theorem \ref{theorem3} we get
$$\lim_{t\rightarrow +\infty} \|\theta(t)\|_{H^s}=0.$$
\hfill$\blacksquare$
\section{\bf Proof of Corollary \ref{Corollary2.5}}
By Theorem \ref{theorem1.4}, their exist a unique global solution $\theta$ of \eqref{AQG} 
 such that $\theta\in  C((0,+\infty),H^2(\R^2)).$ For that, let $t_0>0$ and considering now  the following system
	\begin{equation}\label{AQG2}\tag{$AQG_2$}
		\begin{cases}
			\partial_t\gamma+ u_\gamma.\nabla\gamma +|\partial_1|^{2\alpha}\gamma+ |\partial_2|^{2\beta}\gamma=0,\\
			\gamma(0)=\theta(t_0)\in H^2(\R^2).
		\end{cases}
	\end{equation}
	Therefore, by Theorem \ref{Theorem1.1}, there exists a unique global solution $\gamma$ of \eqref{AQG2} satisfy:
	$$\gamma\in C_b(\R^+,H^2(\R^2)),\quad |\partial_{1}|^\alpha\gamma,|\partial_{2}|^\beta\gamma\in L^2(\R^+,H^2(\R^2)).$$
	By the uniqueness of solution we have 
	$$\theta(t)=\gamma(t-t_0),\quad \forall t\geq t_0.$$
So Using Theorem \ref{theorem3} and Theorem \ref{theorem4}, we obtain
\begin{equation*}
	\lim_{t\rightarrow +\infty}\|\theta(t)\|_{H^2}=0.
\end{equation*}
The fact that 
\begin{equation*}
	\|\theta(t)\|_{H^s}\leq \|\theta(t)\|_{H^2},\quad\forall t\geq t_0.
\end{equation*}
That implies 
\begin{equation*}
	\lim_{t\rightarrow +\infty}\|\theta(t)\|_{H^s}=0.
\end{equation*}
\hfill $\blacksquare$

\end{document}